\documentstyle[11pt,a4]{amsart}

\newtheorem{th}{Theorem}
\newtheorem{lem}{Lemma}

\newtheorem{rem}{Remark}

\def\eps{\varepsilon}

\def\eps{\varepsilon}

\newtheorem{cor}{Corollary}

\title[Regularization for fractional integral...]
{Regularization for fractional integral. Application to nonlinear
equations with singularities
 }

\author[ M. Stojanovi\'{c}]
{Mirjana Stojanovi\'{c}
}

\begin{document}

\maketitle

\footnotetext[1]{
 Department of Mathematics and Informatics, Faculty of science, University of Novi Sad,\\
 Trg D.Obradovi\'{c}a 4, 21 000 Novi Sad, Serbia and Montenegro\\
stojanovic@@im.ns.ac.yu
 }

\footnotetext[2]{ Subject Classificassion (2000): 46F30, 26A33,
45D05, 35K55, 35Q55.

 }

\begin{abstract}
We give the regularization for fractional integral by delta
sequence and apply it to obtain existence-uniqueness theorems in
Colombeau algebras
 for
 nonlinear equations with singularities:
nonlinear system of integral equations with polar kernel and
nonlinear parabolic equations (of ordinary type, with nonlinear
conservative term and
 with Schr\"odinger kernel) with strongly singular initial data
and non-Lipschitz nonlinearities. In a case of nonlinear parabolic
equations we do in fact regularization of heat semigroup with
delta sequence with respect to the time variable $t.$ We do the
same for linear Schr\"odinger equation.

\end{abstract}

\section{\bf Introduction}

Fractional time derivatives have been considered in many problems
in mechanics with application to ordinary fractional differential
equations in mechanics (cf. \cite{hamiltonian}).

In recent years  considerable interest to fractional calculus has
been stimulated by various applications  of that calculus in
numerical analysis and different areas of physics and engineering
possibly including fractal phenomena (cf. \cite{mainardi}).

The purpose of this paper is study of effects of fractional
derivative term in system of  nonlinear integral equations,
nonlinear parabolic equations and linear Schr\"odinger equation.

We find a fractional derivative term in system of nonlinear
integral equations with polar kernel and in heat kernel in
nonlinear parabolic equations as well as in Schr\"odinger
semigroup in linear Schr\"odinger equation. We regularize it by
delta sequence with respect to $t$
 to obtain
existence-uniqueness theorems in Colombeau vector type spaces.

The content of the paper is the following.

 We give the regularization for fractional integral by delta
sequence and obtain the logarithmic boundedeness of fractional
derivative of delta distribution. Then, we apply it to the system
of nonlinear Volterra type integral equations with polar kernel
and singularities for free term to obtain the existence-uniqueness
theorems in Colombeau space $({\cal G}(I))^n,$ where $I$ is an
interval around zero.

We  give the regularization for heat kernel by delta sequence with
respect to the time variable $t$ to nonlinear parabolic equations
of ordinary type, parabolic equations with nonlinear conservative
term and nonlinear parabolic equations with Schr\"odinger kernel
with strongly singular initial data  and non-Lipschitz
nonlinearities.  We apply the results for regularized fractional
integral
 to obtain the existence-uniqueness theorems in
Colombeau vector type spaces ${\cal G}_{C^1,(L^p,L^q)}([0,T), {\bf
R}^n),$ $1\leq p,q \leq \infty,$
 for those choices of $p,q$
and the dimension space $n$
  for which corresponding Colombeau
 space is an algebra with multiplication.
 For example, it holds for  $p=q=2$, $n \leq 3$ due to Sobolev imbedding theorems.

We do the same for linear Schr\"odinger equation with strongly
singular potential and the initial data.

Let us mention that the pioneering work concerning the
regularization of semigroups of nonlinear PDEs
 with respect to the
time variable $t$ in Colombeau algebra was done in
\cite{biagioni+pilipovic}.

\section{\bf Colombau algebras}

For general theory of Colombeau generalized functions and for the
algebra ${\cal G}(\Omega)$, where $\Omega$ is an open set,
 cf. \cite{colombeau1}, \cite{colombeau2},
 \cite{biagioni}, \cite{majklbuk} and recently
\cite{grosser}, \cite{vickers}. We recall the basic definitions
from \cite{grosser}.

Suppose that $\Omega$ is an open set.  Define ${\cal
E}(\Omega)=(C^{\infty}(\Omega))^I,$
$${\cal E}_M(\Omega)=\{(u_{\eps})_{\eps} \in {\cal E}(\Omega)|
\forall K \subset \subset \Omega \forall \alpha \in {\bf N}^n_0
\exists N\in {\bf N}\; \mbox{with}\\ \sup_{x\in K}
|\partial^{\alpha}u_{\eps}(x)|=O(\eps^{-N})$$ as $\eps \to 0 \},$
$$
{\cal N}(\Omega)=\{(u_{\eps})_{\eps} \in {\cal E}(\Omega)| \forall
K \subset \subset \Omega \forall \alpha \in {\bf N}_0^n \forall s
\in {\bf N}:\\
\sup_{x\in K} |\partial^{\alpha}u_{\eps}(x)|=O(\eps^s)\; \mbox{as}
\; \eps \to 0\}.$$ Elements of ${\cal E}_M(\Omega)$ and ${\cal
N}(\Omega)$ are called moderate, resp. negligible functions. The
Colombeau algebra on $\Omega$ is defined  as $ {\cal
G}(\Omega)={\cal E}_M(\Omega)/{\cal N}(\Omega).$

 Recall the construction of the ${\cal G}_{p,q}(\Omega)$ algebras
from \cite{biagioni+ober+siam}.

Let $\Omega \subset {\bf R}^n$ be open, $m \in {\bf Z}$, $1 \leq p
\leq \infty.$ Denote by $W^{m,p}(\Omega)$ usual Sobolev space
whose all derivatives up  to the order $m$ are finite in
corresponding norm. Then,
$$W^{\infty,p}(\Omega)=\cap_{m} W^{m,p}(\Omega), \;
W^{-\infty, p}(\Omega)=\cup_{m} W^{-m, p}(\Omega).$$ Let $1 \leq
p,q \leq \infty.$ Define
$${\cal E}(\Omega)=\{ u| (0,\infty) \times \Omega \to {\bf R} \;
\mbox{such that} \; u_{\eps}(x) \; \mbox{is} \; C^{\infty} \;
\mbox{in} \; x\in \Omega, \; \forall \eps > 0 \}$$
$$
{\cal E}_p(\Omega)=\{ u \in {\cal E}(\Omega), \; \mbox{such that}
\; u_{\eps}  \in W^{\infty, p}(\Omega), \; \forall \eps >0 \}$$
$$
{\cal E}_{M,p}(\Omega)=\{ u \in {\cal E}_p, \; \mbox{such that} \;
\forall \alpha \in {\bf N}^n_0 \exists N \in {\bf N} \; \mbox{such
that} \; ||\partial^{\alpha}u_{\eps}(\cdot)||_p=O(\eps^{-N}), \eps
\to 0 \}$$
$$
{\cal N}_{p,q}(\Omega)=\{ u \in {\cal E}_{M,p}(\Omega) \cap {\cal
E}_q (\Omega) \; \mbox{such that} \; \forall \alpha \in {\bf
N}^n_0 \forall  s \in {\bf N}, \;
||\partial^{\alpha}u_{\eps}(\cdot)||_q=O(\eps^s),\; \eps \to 0
\}$$ where $||\cdot ||_p$ denotes $L^p$-norm. For properties of
these spaces see \cite{biagioni+ober+siam}. Colombeau space
${\cal G}_{p,q}([0,T), {\bf R}^n)$, $1\leq p,q \leq \infty$ is the
factor set
$$
{\cal G}_{p,q}([0,T), {\bf R}^n)={\cal E}_{M,p}([0,T), {\bf
R}^n)/{\cal N}_{p,q}([0,T), {\bf R}^n).$$

We shall give the construction of Colombeau  vector type spaces
${\cal G}_{C^1, (L^p,L^q)}([0,T), {\bf R}^n)$, $1\leq p,q \leq
\infty.$ All of these spaces are not algebra with multiplication.
For $n\leq 3,$ $p=q=2$ due to the Sobolev imbedding $L^2({\bf
R}^n) \subset L^{\infty}({\bf R}^n),$ the Colombeau vector type
space ${\cal G}_{C^1, (L^2,L^2)}([0,T), {\bf R}^n)$ is an algebra
with multiplication. For the proof cf. \cite{n+p+r}.

 Let $\Omega$ be an open set in
${\bf R}^n$, $1\leq p,q \leq \infty.$

Space ${\cal E}_{C^1,(L^p,L^q)}([0,T), {\bf R}^n)$ is a vector
space of nets $(G_{\eps})_{\eps}$ such that
$$
G_{\eps}\in C^0([0,T), L^p({\bf R}^n)) \cap C^1((0,T), L^p({\bf
R}^n)), \; \eps <1,$$ with the following property: $\forall T_1
\in (0,T) \exists N \in {\bf N}$ such that
\begin{equation}
\label{casa}
 \max\{ \sup_{t\in [0,T)} ||G_{\eps}(t)||_{L^{p}},
\; \sup_{t\in [T_1,T]} ||\partial_t G_{\eps}(t)||_{L^p} \}
=O(\eps^{-N}), \; \eps \to 0.
\end{equation}
If
$$
 \max\{ \sup_{t\in [0,T)} ||G_{\eps}(t)||_{L^{q}},
\; \sup_{t\in [T_1,T]} ||\partial_t G_{\eps}(t)||_{L^q} \}
=O(\eps^s), \; \eps \to 0,
$$
 holds for $\forall s \in {\bf N}$ we obtain null space
 ${\cal
N}_{C^1, (L^{p},L^{q})}([0,T), {\bf R}^n).$ Colombeau vector type
space is defined as  the quotient space
$$
{\cal G}_{C^1, (L^{p}, L^q)}([0,T), {\bf R}^n)={\cal
E}_{C^1,L^{p}}([0,T), {\bf R}^n)/{\cal N}_{C^1, (L^{p},
L^q)}([0,T), {\bf R}^n).$$

\section{\bf Basic facts from fractional calculus}

Fractional integrals are defined as a generalization of Cauchy
formula for a repeated indefinite integral
$$
J^nf(t):=f_n(t)=\frac{1}{(n-1)!} \int_0^t (t-\tau)^{n-1}
f(\tau)d\tau, \; t>, \; n\in {\bf N},$$ where $Jf$ is a primitive
function of $f:$
$$(Jf)(t):=\int_0^t f(\tau)d\tau,$$
and $J^n$ denotes the $n^{th}$ power of $J.$

For an arbitrary real number $\alpha >0$ and sufficiently regular
function $f$ on ${\bf R}$ we define fractional integral of order
$\alpha >0$
\begin{equation}
\label{sesnaest} J^{\alpha}f(t):=\frac{1}{\Gamma(\alpha)} \int_0^t
(t-\tau)^{\alpha-1} f(\tau)d\tau, \; t>0, \; \alpha \in {\bf R}^+.
\end{equation}
The equation (\ref{sesnaest}) can be expressed in terms of
distributions for arbitrary complex number $\alpha$ (cf.
\cite{geljfandshilov}),
 as
\begin{equation}\label{sterio}
{\bf \Phi}_{\alpha}(t):= \left\{\begin{array}{ll}
\frac{t_+^{\alpha-1}}{\Gamma(\alpha)}, \; \; \alpha >0,
\\
D^{(n)}{\bf \Phi}_{\alpha+n}, \; \alpha \leq 0, \; \alpha+n>0, \;
n\in {\bf N},\end{array}\right.
\end{equation}
(where $D^{(n)}$is the $n^{th}$ distributional derivative) in the
form $ J^{\alpha}f:={\bf \Phi}_{\alpha}*f_+,$ where
$f_+(t):=H(t)f(t)$, $H(t)$ is the Heaviside function, convolution
is defined in distributional sense. The distribution
$t_+^{\alpha-1}$ for $\alpha >0$ is defined by
$$
t_+^{\alpha-1}=\left\{\begin{array}{ll} 0 \;\; \; \; \; \; t<0\\
t^{\alpha-1} \;\; \; \;  t>0.\end{array}\right.$$ For other values
of $\alpha$ this distribution is defined by analytic continuation
or by regularization of divergent integrals (cf.
\cite{geljfandshilov}). We have
$$
J^{\alpha}f(t)={\bf \Phi}_{\alpha}*f(t), \; \alpha >0, \;
\mbox{and} \; D{\bf \Phi}_{\alpha}={\Phi}_{\alpha-1},$$ and
semigroup property holds
$${\bf \Phi}_{\alpha}(t)*{\bf \Phi}_{\beta}(t)={\bf
\Phi}_{\alpha+\beta}(t), \; \alpha, \beta >0.$$

With  interpretation of the quotient as a limit if $t=0$ (cf.
\cite{mainardi})
$$
{\bf \Phi}_{-n}:=\frac{t_+^{n-1}}{\Gamma(-n)}=\delta^{(n)}(t), \;
n\in {\bf N}_0,$$ where $\delta^{(n)}(t)$ is the generalized
derivative of order $n$ of the Dirac delta distribution. Then,
$$
\frac{d^n}{dt^n}f(t)=f^{(n)}(t)={\bf
\Phi}_{-n}(t)*f(t)=\int_{0^-}^{t^+}
f(\tau)\delta^{(n)}(t-\tau)d\tau, \; t>0,$$ since
$$
\int_{0^-}^{t^+}f(\tau) \delta^{(n)}(\tau-t)d\tau=(-1)^n f^n(t),
\; \delta^{(n)}(t-\tau)=(-1)^n \delta^{(n)}(\tau-t).$$ The formal
definition of fractional derivative could be
\begin{equation}\label{zvezda}
{\bf \Phi}_{-\alpha}(t)*f(t)=\frac{1}{\Gamma(-\alpha)}
\int_{0^-}^{t^+}\frac{f(\tau)}{(t-\tau)^{\alpha+1}} d\tau, \;
\alpha \in {\bf R}^+.
\end{equation}
In general, ${\bf \Phi}_{-\alpha}(t)$ is not locally absolutely
integrable and integral is divergent.

In order to obtain the definition that is still valid for
classical functions we regularize the divergent integral by delta
sequence. Instead of the function $f$ in (\ref{zvezda}) we shall
use the mollifier $\phi_{\eps}(t)=|ln\eps| \phi(x \cdot |ln\eps|)$
(resp. using $(ln|ln\eps|))$ where $\phi(t)\in C_0^{\infty}({\bf
R}),$
 $\phi(t)\geq 0,$ $\int \phi(t)dt=1.$

 We shall consider three different cases with various values of
 $\alpha.$

 {\bf 1.}
 Consider (\ref{zvezda}).
 For the elements of the set $n \in {\bf N}_0$,
 ${\bf
 \Phi}_{-n}(t)*\phi_{\eps}(t)=\delta^{(n)}(t)*\phi_{\eps}(t).$
 We shall find $L^1$-norm of the above convolution for later use.
 We have
 $$
 ||{\bf
 \Phi}_{-n}(t)*\phi_{\eps}(t)||_{L^1}=||\delta^{(n)}(t)*\phi_{\eps}(t)||_{L^1}$$
 $$
 =||\delta(t)*D^{(n)}\phi_{\eps}(t)||_{L^1}
=
 ||D^{(n)}\phi_{\eps}(t)||_{L^1}\leq C |ln\eps|^n$$
 where $D^{(n)}$ denotes distributional derivative.

 {\bf 2.}
 Let $\alpha >0.$
 Then,
 ${\bf \Phi}_{\alpha}(t):=\frac{t_+^{\alpha-1}}{\Gamma(\alpha)},$
 $\alpha >0,$
 and fractional integral is given by (\ref{sesnaest}).
 Then,
 $
 J^{\alpha}f(t)={\bf \Phi}_{\alpha}(t)*f(t), \; \alpha>0.$
Then,
 $
 ||J^{\alpha}\phi_{\eps}(\cdot)||_{L^1}\leq ||{\bf
 \Phi}_{\alpha}(\cdot)||_{L^1}
 ||\phi_{\eps}(\cdot)||_{L^1}\leq C.$
 By this, we cover the case $0< \alpha <1.$

 {\bf 3.}
 Consider the case $\alpha <0.$
 By classical definition
(\ref{sterio})
 of the function ${\bf
 \Phi}_{\alpha}(t)$
for $\bar{\alpha} \in {\bf R}_+,$ $\alpha=-\bar{\alpha},$
$\bar{\alpha}>0,$
 we have  ${\bf
\Phi}_{-\bar{\alpha}}(t)=D^{(n)}{\bf \Phi}_{-\bar{\alpha}+n}(t),$
$\bar{\alpha} \geq 0,$ $-\bar{\alpha}+n>0,$ $n\in {\bf N},$
$(\bar{\alpha} <n).$ Then,
$$J^{-\bar{\alpha}}\phi_{\eps}(t)=\frac{1}{\Gamma(-\bar{\alpha})}\int_{0^-}^{t^+}
(t-\tau)^{-\bar{\alpha}-1}\phi_{\eps}(\tau)d\tau, \; \bar{\alpha}
\in {\bf R}^+, \; t>0,$$ and
$$
J^{-\bar{\alpha}}\phi_{\eps}(t)={\bf
\Phi}_{-\bar{\alpha}}*\phi_{\eps}(t)=D^n{\bf
\Phi}_{-\bar{\alpha}+n}(t)*\phi_{\eps}(t)= {\bf
\Phi}_{-\bar{\alpha}+n}(t)*D^n\phi_{\eps}(t), \;
-\bar{\alpha}+n>0.$$ By case ${\bf 2.}$
$$
||J^{-\bar{\alpha}}\phi_{\eps}(t)||_{L^1}\leq ||{\bf
\Phi}_{-\bar{\alpha}+n}(t)||_{L^1} ||D^n\phi_{\eps}(t)||_{L^1}\leq
C |ln\eps|^n.$$ Thus, we have for each $\alpha \in {\bf R}$
\begin{equation}
\label{trizvezdice} |||x|_+^{\alpha-1}*\phi_{\eps}(x)||_{L^1}\leq
 \left\{\begin{array}{ll} C \;
\;\; \mbox{when} \; \alpha >0,\\
C|ln\eps|^m \; \; \mbox{when}\; \; \alpha \leq 0, \;
m>\bar{\alpha}, \; \alpha=-\bar{\alpha}, \; \bar{\alpha}>0.
\end{array}\right.
\end{equation}
We have just proved the following Lemma.
\begin{lem}\label{lemmaone}Fractional integral
of the delta sequence $$
J^{\alpha}\phi_{\eps}(t):=\frac{1}{\Gamma(\alpha)} \int_0^t
(t-\tau)^{\alpha-1} \phi_{\eps}(\tau)d\tau, \; t>0, \; \alpha \in
{\bf R},
$$
where $\phi(x)\in C_0^{\infty}({\bf R}),$ $\phi(x)\geq 0,$ $\int
\phi(x)dx=1,$ and $\phi_{\eps}(x)=|ln\eps|\phi(x \cdot |ln\eps|)$
(resp. using $(ln|ln\eps|)$ instead of $|ln\eps|)$ has  the
following bounds in $L^1$-norm:
$$||J^{\alpha}(\phi_{\eps}(t))||_{L^1}=
|||t|_+^{\alpha-1}*\phi_{\eps}(t)||_{L^1}\leq
 \left\{\begin{array}{lll} C \;
\;\; \mbox{when} \; \alpha >0,\\
C|ln\eps|^m \; \; (\mbox{resp.} \; \; (ln|ln\eps|)^m)\; \;
\mbox{when}\\ \alpha \leq 0, \; m>\bar{\alpha}, \; \;
\alpha=-\bar{\alpha}, \; \bar{\alpha}>0.
\end{array}\right.$$
\end{lem}
\begin{cor}
When $\alpha=0$ we have $||vp
\frac{1}{|x|_+}*\phi_{\eps}(x)||_{L^1}\leq C |ln\eps|^m, \; m>0.$
\end{cor}

\section{\bf Application to nonlinear Volterra integral equation}

Consider a system of nonlinear Volterra integral equations with
polar kernel
\begin{equation}\label{laz}
 f^i(x)=g^i(x)+\int_0^x
\frac{K^i(x,y,f(y))}{|x-y|^{-\alpha+1}}dy, \; x\in I, \;
i=1,..,n,\; \alpha \in {\bf R},
\end{equation}
 where $I$ is an interval around
zero,
 $g^i(x) \in {\cal D}'(I),$ $K^i \in L_{loc}^{\infty}(I\times I \times {\bf R}^{2n})$
  and it is
not of  Lipschitz class, $i=1, 2, ... n$.
\begin{th}
Nonlinear system of integral equations (\ref{laz}) with three type
 of singularities: $g(x)\in ({\cal D}'(I))^n$, $ K(x,y,f(y)) \in
 (L_{loc}^{\infty}(I\times I\times {\bf R}^{2n}))^n$
 is nonlinear and non-Lipschitz,
 and integral is divergent on the set $\{ y|y=x\}$
 has a unique
solution in the space $[f_{\eps}(x)]\in ({\cal G}(I))^n.$
\end{th}
{\bf Proof.}
The function $K^i$ is substituted by $K_{\eps}^i$ and $g^i$ by
$g_{\eps}^i$, $i=1,...,n$, $\eps \in (0,1]$, in order to avoid
non-Lipschitz nonlinearity of $K^i$ and singularities of $g^i$,
$i=1,...,n.$

For given $g=(g^1,...,g^n)\in ({\cal D}'(I))^n,$
we put $g_{\eps}^i=(g^i \kappa_{\eps})*\phi_{\eps},$ $i=1,...,n,$
where $\kappa_{\eps}\in C_0^{\infty}(I),$
$$\kappa_{\eps}=\left\{\begin{array}{ll} 1 \; \mbox{on} \;
I_{2\eps}\\
0 \; \mbox{on} \; I\setminus I_{\eps},
\end{array}\right.$$
$I=(-a,a)$, $I_{j\eps}=(-a+j\eps, a-j\eps),$ $j=1,2$ and
$\phi_{\eps}(x)=|ln\eps|\phi(x\cdot |ln\eps|),$ $x\in {\bf R},$ is
mollifier with the properties: $\phi \in {\cal S}({\bf R}),$ $\int
\phi(x)dx=1,$ $\int x^{\beta}\phi(x)dx=0,$ for  all $\beta \in
{\bf N}_0^n.$ If $a=-\infty,$ or $b=\infty$ then
$a+j\eps=-\infty,$ or $b-j\eps=\infty.$

By cut-off method (cf. \cite{polar}) we have
$$
|\nabla K_{\eps}(x,y,f_{\eps}(y))| \leq C |ln\eps|^b, \;
(\mbox{resp.}\;
 \leq C (ln|ln\eps|)),$$ where $b$ will be determined to
handle the problem under consideration. Then, the integral is
singular only due to $|x-y|^{\alpha}.$

Consider the convolution form of the above equation
\begin{equation}\label{cetirizvezdice}
f_{\eps}(x)=g_{\eps}(x)+
K_{\eps}(x,y,f_{\eps}(y))*|x|_+^{\alpha-1}*\phi_{\eps}(x),
\end{equation}
where we regularized fractional part of the integrand by delta
sequence to avoid the divergence of the integrand along the
diagonal $x=y.$ In integral form we have
\begin{equation}\label{sedam}
f_{\eps}(x)=g_{\eps}(x)+\int_0^x \int_{{\bf R}}
\frac{K_{\eps}(x,y,f_{\eps}(y))}{|x-y|^{-\alpha+1}}
\phi_{\eps}(x-y-s)dsdy.\end{equation}
 By H\"older inequality
$$
|f_{\eps}(x)| \leq C|ln\eps|^{1}+\int_0^x
||K_{\eps}(x,y,f_{\eps}(y)||_{L^{\infty}}
||\phi_{\eps}(x-y)*|x-y|_+^{\alpha-1}||_{L^1}dy.$$ Let $\alpha\leq
0.$ By (\ref{trizvezdice})
$$
f_{\eps}(x)=g_{\eps}(x)+\int_0^x ||\nabla K_{\eps}(x,y,\theta
f_{\eps}(y))||_{L^\infty} \sup_y |f_{\eps}(y)| |ln\eps|^m dy, \;
-\bar{\alpha}+m>0,$$ where $\alpha=-\bar{\alpha}, \bar{\alpha}\geq
0.$
Gronwall inequality yields
$$
|f_{\eps}(x)|\leq C |ln\eps|^{1} \exp{(C|ln\eps|^{b+m})}\leq C
\eps^{-N}, \; \exists N>0,$$ where $m \geq \bar{\alpha},$ i.e.
$b+\bar{\alpha} <1.$ For $\alpha >0$ due to (\ref{trizvezdice}) we
have only the condition $b<1.$

Consider $\beta^{th}$-derivative, $\beta \in {\bf N}_0,$ $\beta
\geq 1.$
 We have
in (\ref{sedam})
$$D^{\beta}f_{\eps}(x)=D^{\beta}g_{\eps}(x)+\int_0^x \int_{{\bf
R}} K_{\eps}(x,y,f_{\eps}(y))
\frac{\partial_x^{\beta}\phi_{\eps}(x-y-s)}{|x-y|^{-\alpha+1}}dsdy.$$
By H\"older inequality we obtain
\begin{equation}\label{pap}
D^{\beta}f_{\eps}(x)=D^{\beta}g_{\eps}(x)+\int_0^t \sup_{y}
|\nabla K_{\eps}| \sup_{y} |f_{\eps}(y)|
||\partial_x^{\beta}\phi_{\eps}(x-y)*|x-y|_+^{\alpha-1}||_{L^1}
dy.\end{equation}
 We
shall calculate the first
$$
||\partial_x^{\beta}\phi_{\eps}(x)*|x|_{+}^{\alpha-1}||_{L^1} \leq
\int_{\bf
R}|(\phi_{\eps}(x)*\partial_x^{\beta}|x|_+^{\alpha-1})|dx$$ $$=
\int_{\bf R} \int_{\bf R}
|(\alpha-1)(\alpha-2)...(\alpha-\beta-1)| |\phi_{\eps}(s)|
|x-s|_{+}^{\alpha-1-\beta}dsdx.$$ By Lemma \ref{lemmaone}
 we have
 \begin{equation}\label{australia}
 ||\phi_{\eps}(x)*\partial_x^{\beta}|x|_+^{\alpha-1}||_{L^1} \leq
 \left\{\begin{array}{lll}
 C \; \; \mbox{when} \; \;\alpha >0, \; \alpha>\beta,\\
C|ln\eps|^m, \; \; \alpha >0, \; \alpha<\beta, \; m>-\beta+\alpha
  \\
 C |ln\eps|^m \; \; \alpha<0, \; \alpha-\beta<0, \; \;  m>-\bar{\alpha}-\beta,
  \; \alpha=-\bar{\alpha},
 \bar{\alpha}>0.\end{array}\right.
 \end{equation}
Setting this in (\ref{pap}) we obtain for $\alpha \leq 0$
$$
D^{\beta}f_{\eps}(x)= D^{\beta}g_{\eps}(x)+\int_0^x |ln\eps|^b
\sup_{y} |f_{\eps}(y)| |ln\eps|^m dy, \; m>-\bar{\alpha}-\beta, \;
\bar{\alpha}<\beta.$$ Since $g(x)\in {\cal D}'(I)$ we shall give a
regularization as in (\cite{polar})
$||\kappa_{\eps}f(x)*\phi_{\eps}(x)||_{L^{\infty}} \leq C
|ln\eps|.$ Then Gronwall inequality yields
$$
|D^{\beta}f_{\eps}(x)| \leq C |ln\eps|^{\beta+1} +\int_0^x
|ln\eps|^{b+m} \sup_y |f_{\eps}(y)| dy, \; m> -\bar{\alpha}-\beta,
\; \bar{\alpha} < \beta, \; \bar{\alpha}>0,$$ and
$$
|D^{\beta}f_{\eps}(x)| \leq C |ln\eps|^{1+\beta}
\exp{(C|ln\eps|^{b+m})}\leq C\eps^{-N}, \; \exists N>0, \; x\in
{\bf R}, \; \eps<\eps_0,$$ where $m>-\bar{\alpha}-\beta, $
$b+m<1,$ $b<1+\bar{\alpha}+\beta.$ This condition is included  in
the condition imposed at the first step $m>\bar{\alpha},$
$b<1-\bar{\alpha}.$

If $\alpha >0$ we have by (\ref{australia})
$$
D^{\beta}f_{\eps}(x)=D^{\beta}g_{\eps}(x)+\int_0^x C
\left\{\begin{array}{ll} C \; \; \alpha>0, \; \alpha>\beta\\
C |ln\eps|^m, \; \; \alpha>0, \; \alpha < \beta, \;
m>-\beta+\alpha\end{array}\right. |ln\eps|^b \sup_{y}
|f_{\eps}(y)|dy.$$ Gronwall inequality yields
$$
|D^{\beta}f_{\eps}(x)| \leq C |ln\eps|^{1+\beta} \exp{(C\left\{\begin{array}{ll}
C \; \; \alpha>0, \; \alpha>\beta\\
C |ln\eps|^m, \; \; \alpha>0, \; \alpha < \beta, \;
m>-\beta+\alpha\end{array}\right. |ln\eps|^b)}\leq C  \eps^{-N},$$
$ \exists N>0, \; x\in {\bf R}, \; \eps <\eps_0,$ under the
condition $b<1,$ when $\alpha>\beta,$ and $b+m<1,$
$b<1+\beta-\alpha$, $\alpha<\beta,$ what is included in the first
step condition $b<1.$
 Thus, we have the moderateness of the solution
$f_{\eps}(x) \in {\cal E}_M({I}).$

Let us prove the uniqueness. Suppose  that $f_{1\eps}(x)$ and
$f_{2\eps}(x)$ are two solutions to the equation (\ref{sedam}) and
denote their difference by $F_{\eps}(x).$ Then, we have
\begin{equation}\label{jelen}
|F_{\eps}(x)| \leq \int_0^x
|K_{\eps}(x,y,f_{1\eps}(y))-K_{\eps}(x,y,f_{2\eps}(y))|
\frac{|\phi_{\eps}(x-y-s)|}{|x-y|_+^{1-\alpha}}ds dy+d_{\eps}(x),
\end{equation}
where $d_{\eps}(x)\in {\cal N}(I). $ By mean value theorem
$$
|K_{\eps}(x,y,f_{1\eps}(y))-K_{\eps}(x,y,f_{2\eps}(y))|\leq C
|f_{1\eps}(y)-f_{2\eps}(y)|$$ $$\int_0^1 |\nabla
K_{\eps}(x,y,\theta f_{1\eps}(y)+(1-\theta)f_{2\eps}(y))| d\theta
\leq C |F_{\eps}(y)| |ln\eps|^b.$$ Setting this in (\ref{jelen})
we obtain
$$
|F_{\eps}(y)|\leq C \int_0^x \sup_y |F_{\eps}(y)| |ln\eps|^b
||\phi_{\eps}(x-y)*|x-y|_+^{\alpha-1}||_{L^1} dy+d_{\eps}(x).$$
Applying Lemma \ref{lemmaone}
 we obtain for $\alpha=-\bar{\alpha}<0,$
 where $\bar{\alpha}>0$
 $$
 |F_{\eps}(x)| \leq C \int_0^x \sup_y |F_{\eps}(y)| |ln\eps|^{b+m}
 dy+d_{\eps}(x), \; m>\bar{\alpha}.$$
Gronwall inequality yields
$$
|F_{\eps}(x)| \leq C d_{\eps}(x) \exp{(C |ln\eps|^{b+m})}\leq C
\eps^s, \; \forall s \in {\bf N}, \; x\in {\bf R}, \; \eps<1,$$
when $b+m<1,$ $m>\bar{\alpha},$ i.e. $b+\bar{\alpha}<1.$

If $\alpha>0$ by Lemma \ref{lemmaone} we have
$$
|F_{\eps}(x)| \leq C \int_0^x \sup_y |F_{\eps}(y)| |ln\eps|^b
dy+d_{\eps}(x).$$ Gronwall inequality yields
$$
|F_{\eps}(x)| \leq C d_{\eps}(x) \exp{(C |ln\eps|^b)}\leq C
\eps^s, \; \forall s \in {\bf N}, \; x\in {\bf R}, \;
\eps<\eps_0,$$ when $b<1.$

Consider $\beta^{th}$-derivative, $\beta \in {\bf N}_0,$ $\beta
\geq 1,$
 applied
on $F_{\eps}(x)$ in the equation (\ref{pap})
to obtain
$$|D^{\beta}F_{\eps}(x)| \leq \int_0^x |\int_0^1 \nabla
K_{\eps}(x,y, \theta f_{1\eps}+(1-\theta)f_{2\eps})d\theta|
$$
$$\sup_y |F_{\eps}(y)|
(||\phi_{\eps}(x-y)*\partial_x^{\beta}|x-y|_+^{\alpha-1}||_{L^1})dy+d_{\eps}(x),
\; 0<\theta <1.$$ By (\ref{australia}) we have
$$
|F_{\eps}^{\beta}(x)| \leq C \int_0^x |ln\eps|^b \sup_y
|F_{\eps}(y)| \left\{\begin{array}{lll} C \; \; \mbox{when} \;
\alpha > \beta, \; \alpha >0 \\ C |ln\eps|^m, \; \mbox{when} \;
\alpha>0, \; \alpha-\beta <0, \;  m
>\alpha-\beta
\\
C|ln\eps|^m, \; \alpha<0, \; \alpha-\beta<0, \;
m>-\bar{\alpha}-\beta
 \end{array}\right.$$
 $$
dy+d_{\eps}(x).$$ Using the null properties of $F_{\eps}(x)$ from
the first step of the induction we obtain
$$
|D^{\beta}F_{\eps}(x)| \leq C d_{\eps}(x)+
(\left\{\begin{array}{lll} C \; \; \mbox{when} \; \alpha > \beta
\; \alpha >0 \\ C |ln\eps|^m, \; \mbox{when} \; \alpha>0 \;
\alpha-\beta <0, \;  m
>\alpha-\beta
\\
C|ln\eps|^m, \; \alpha<0, \; \alpha-\beta<0, \;
m>-\bar{\alpha}-\beta
 \end{array}\right.)\leq C \eps^s, $$ $
\forall s \in {\bf N}, \; \eps<\eps_0, \; x\in {\bf R},$ under the
condition
$$\left\{\begin{array}{lll}
b<1, \; \;\alpha-\beta>0, \; \alpha>0\\
b+m<1, \; m>\alpha-\beta, \; \alpha>0, \; \alpha-\beta<0\\
b+m, \; m>-\bar{\alpha}-\beta, \; \alpha<0, \;
\alpha-\beta<0\end{array}\right.$$ which is  included in the
condition  imposed at the first step
$$
\left\{\begin{array}{ll} b<1, \; \alpha>0\\
b<1-\bar{\alpha}, \; \; \alpha<0.
\end{array}\right.$$
 Then, $F_{\eps}(x)\in {\cal N}({I}).$
 Follows $f_{1\eps}(x) \approx f_{2\eps}(x),$
 the solution is unique.

 Thus, there exists a unique solution in the space $[f_{\eps}] \in
 ({\cal G}({I}))^n.$
$\Box$

\begin{rem}
 The condition imposed on parameters gives
a link between the singularity of the fractional integral and
non-Lipschitz nonlinearity of the kernel $K.$ The singularity of
the fractional integral is reflected by the growth of the
mollifiers. The singularity of the free term is reflected by
$(|ln\eps|)$-boundedeness which is necessary for moderateness in
  applying Gronwall
inequality. When we consider $(ln|ln\eps|)$-boundedeness for both
$K$ and fractional part, the slower growth of mollifiers
contributes to the moderateness.
\end{rem}

\section{\bf Application to nonlinear parabolic equations with
strongly singular initial data}\label{subscript}

Consider nonlinear parabolic equations with strongly singular
initial data and non-Lipschitz nonlinearity for free term for
ordinary nonlinear parabolic equations (cf.
\cite{biagioni+gramchev1}), equations with nonlinear conservative
term (cf. \cite{biagioni+gramchev2}) and parabolic equation with
Schr\"odinger kernel (cf. \cite{n+p+r}).

We shall prove the existence-uniqueness theorems in the Colombeau
vector type spaces ${\cal G}_{C^1,(L^p,L^q)}([0,T), {\bf R}^n),$
$1\leq p,q \leq \infty,$
 using the regularization for fractional part of the heat
kernel by delta sequence with respect to the time variable $t.$

We shall consider:

{\bf 1.} Cauchy problem for nonlinear parabolic equation
\begin{equation}\label{greska}
\partial_t u =\triangle u + g(u), \; t>0, \; x\in {\bf R}^n, \;
u(0,x)=\mu(x),
\end{equation}
where $g(u)\in L_{loc}^{\infty}([0,T), {\bf R}^n)$ and does not
satisfy Lipschitz condition, $\mu(x)\in {\cal M}^k({\bf R}^n)$ $
\subset $ $ {\cal D}'$ $({\bf R}^n)$, $k \in {\bf Z}_+$, where
${\cal M}^k({\bf R}^n)$ is the strong dual of Banach space
$C_b^k({\bf R}^n)$ of all $C^k({\bf R}^n)$ functions with bounded
derivatives up to the order $k$.
For example,
 strongly singular initial
data are given in the form of regularization by smooth mollifiers
of sum of derivatives of Dirac measures. We put some conditions on
mollifiers regularizing the leading term of the initial data:
$$\phi_{\eps}\in C_0^{\infty}({\bf R}^n), \; \int_{{\bf R}^n}
\phi_{\eps}(x)dx=1, \; \lim_{\eps \to 0}\phi_{\eps}(x)=\delta(x)$$
in ${\cal D}'({\bf R}^n)$ where $\phi_{\eps}(x)=|ln\eps|^{an}
\phi(x \cdot |ln\eps|),$ $\eps>0,$ and $\delta(x)$ stands for the
Dirac measure massed at the point zero (cf. \cite{gramchev}):
$$\lim_{\eps \to 0} u_{0\eps}(x)=u_0(x):=\sum_{i=1}^n \sum_{j=0}^k
a_{ik} \delta^{(j)}(x-\xi_i),\; \mbox{in} \; {\cal D}'({\bf
R}^n).$$

{\bf 2.} Cauchy problem for nonlinear  parabolic equation with
conservative nonlinear term
\begin{equation}\label{studenti}
\partial_t u -\triangle u +\partial_x \cdot \vec{g}(u)=0, \; t>0,
\; x\in {\bf R}^n, \; u(0,x)=u_{0}(x)=D^k\psi(x),
\end{equation}
where $u=(u_1,...,u_m),$ $\vec{g}(u)=(g(u), ..., g_n(u))\in
L_{loc}^{\infty}([0,T), {\bf R}^n)$, $\partial_x \cdot
\vec{g}(u)=\vec{g}(u) \cdot \nabla u=\sum_{j=1}^n g'_j(u)
\partial_{x_j}u,$
$k\geq0,$ where the initial data are the following:
$$
\mu(x)=|D|^k\psi, \; \psi \in L^p({\bf R}^n), \; k>0,$$ $ 1\leq p
\leq \infty, \; D=(-\triangle)^{1/2}, \; \psi \in L^p({\bf R}^n),
\; 1\leq p \leq \infty.$ We use the following regularization for
the initial data:
$$
\mu_{\eps}(x)=D^k\psi(x)*\phi_{\eps}(x)=\psi(x)*D^k\phi_{\eps}(x),
\; \psi \in L^p({\bf R}^n),\;
\phi_{\eps}(x)=|ln\eps|^{an}\phi(x\cdot |ln\eps|),$$ and then
$$||\mu_{\eps}(\cdot)||_{L^p} \leq C
||D^k\phi_{\eps}(\cdot)||_{L^1}\leq C |ln\eps|^{n(a-1)+k}.$$

As the initial data we can consider the powers of delta
distribution $u_0(x) = \delta^{(\beta)}(x),$ $\beta \in
(0,\infty),$ $x\in {\bf R}^n.$
 Without loss of generality
suppose that $u_0(x)=\delta(x),$ $\delta_{\eps}(x)=|ln\eps|^{an}
\phi(x \cdot |ln\eps|),$ $\phi(x) \in C_0^{\infty}({\bf R}^n),$
$\int \phi(x)dx=1,$ $\phi(x)\geq 0.$

{\bf 3.} Nonlinear parabolic equation  with Schr\"odinger kernel,
delta as the potential and delta as the initial data:
\begin{equation} \label{zadaci}
(\partial_t-\triangle)u+V(x)u +g(u)=0, \; u(0,x)=u_0(x), \; x\in
{\bf R}^n,
\end{equation}
where $V(x)$ and $u_0(x)$ are the singular distributions, say
$\delta$-distributions, $g(u)\in L_{loc}^{\infty}([0,T), {\bf
R}^n)$, is nonlinear, non-Lipschitz. The initial data  and
potential could be  sum of derivatives and powers of
$\delta$-distribution. Without loss of generality we shall
consider the case when $u_0(x)=V(x)=\delta(x).$

\subsection{\bf Nonlinear parabolic equation (\ref{greska})}
\label{petjedan}

Consider regularized equation (\ref{greska})
\begin{equation}
\label{jedanaest}
\partial_tu_{\eps}(t,x)=\triangle
u_{\eps}(t,x)+g_{\eps}(u_{\eps}(t,x)), \;
u_{0\eps}(x)=\delta_{\eps}(x),
\end{equation}
where
\begin{equation}\label{dvanaest}
\delta_{\eps}(x)=|ln\eps|^{an} \phi(x \cdot |ln\eps|), \; |\nabla
g_{\eps}(\theta u_{\eps})|\leq C |ln\eps|^b,
\end{equation}
$\phi(x)\in C_0^{\infty}({\bf R}^n),$ $\phi \geq 0,$ $\int
\phi(x)dx=1,$ $g(u)$ is  regularized by cut-off to avoid
non-Lipschitz nonlinearity and the heat semigroup is regularized
as follows.

We shall regularize heat kernel in $L^1$-norm. We have
$$E_{n\eps}(t,x-y)=
E_n(t,x-y)*\phi_{\eps}(t)=\int_{\bf R}
E_n(t-\tau,x-y)\phi_{\eps}(\tau)d\tau$$ where
$$
E_n(t-\tau,x-y)=(4\pi(t-\tau))^{-n/2} \exp{(-|x-y|^2/(4(t-\tau))}.
$$
Then,
$$||
E_{n\eps}(t,x-\cdot)||_{L^1} =\int_{\bf R^n}\int_{\bf R}
|(4\pi(t-\tau))^{-n/2}| |\exp{(-|x-y|^2/(4\pi(t-\tau))}|
|\phi_{\eps}(\tau)|d\tau|dy,$$ and
$$
\partial_x^{\beta}E_{n\eps}(t,x-y)=\partial_x^{\beta}E_n(t,x-y)*\phi_{\eps}(t)=
\int_{\bf R}
\partial_x^{\beta}
E_n(t-\tau,x-y) \phi_{\eps}(\tau)d\tau $$ $$\leq C \int_{\bf R}
(4\pi (t-\tau))^{-n/2} (t-\tau)^{-\beta}
\exp{(-|x-y|^2/(4(t-\tau))}) \phi_{\eps}(\tau)d\tau.$$ Since from
\cite{biagioni+gramchev2} we have
\begin{equation}\label{michael}
||t^{k/2+n/2(1-1/r)} \partial_x^{\beta}E_n(t,\cdot)||_{L^r} \leq
\infty, \; |\beta|\leq k, \; 1\leq r \leq \infty,
\end{equation}
 by Fubini
theorem we obtain
$$||\partial_x^{\beta}E_n(t,x-\cdot)||_{L^1} \leq
C \int_{{\bf R}} \int_{{\bf R}^n} |(t-\tau)^{-n/2-\beta}|
|\phi_{\eps}(\tau)| | \exp{(-|x-y|^2/(4(t-\tau))}|dy d\tau,$$ and
by (\ref{michael}) we continue
$$\leq  C \int_{\bf R} (t-\tau)^{-\beta/2}
||(t-\tau)^{\beta/2}\partial_x^{\beta}E_n(t-\tau,x-\cdot)||_{L^1}
d\tau.$$ Since
$||(t-\tau)^{\beta/2}\partial_x^{\beta}E_n(t-\tau,x-\cdot)||_{L^1}\leq
C$ we obtain
\begin{equation}\label{montreal}
||\partial_x^{\beta}E_{n\eps}(t-\tau, x-\cdot)||_{L^1} \leq
\left\{\begin{array}{ll}
C \; \; \mbox{when} \; \; \beta <2\\
C (ln|ln\eps|)^m, \; \mbox{when} \; \beta \geq 2, \; m> \beta/2-1,
\end{array} \right.
\end{equation}
where we used the $(ln|ln\eps|)$-boundedeness of the heat
semigroup to accomplish the moderateness in corresponding
Colombeau space.

\begin{th}
The regularized equation (\ref{jedanaest}) to  the nonlinear
parabolic equation (\ref{greska}) where $u_0(x)=\delta(x)$, $g(u)
\in L_{loc}^{\infty}([0,T), {\bf R}^n)$ is nonlinear and
non-Lipschitz has a unique solution in Colombeau vector type
spaces ${\cal G}_{C^1, (L^p, L^q)}([0,T), {\bf R}^n)$ for those
choices $1\leq p, q \leq \infty$ for which corresponding Colombeau
space is an algebra with multiplication. For $p=q=2$ the space
${\cal G}_{C^1, (L^2,L^2)}([0,T), {\bf R}^n)$ is an algebra with
multiplication for $n \leq 3$.
\end{th}
{\bf Proof.}
Consider regularized integral form of the nonlinear parabolic
equation (\ref{greska})
$$
u_{\eps}(t,x)=\int_{{\bf R}^n} E_{n\eps}(t,x-y)
u_{0\eps}(y)dy+\int_0^t \int_{{\bf R}^n} E_{n\eps}(t-\tau,x-y)
g_{\eps}(u_{\eps}(\tau,y))dyd\tau,$$ $ t\in [0,T), \; x\in {\bf
R}^n.$
We have in $L^p$-norm
$$
||u_{\eps}(t,\cdot)||_{L^p} \leq ||E_{n\eps}(t,x-\cdot)||_{L^1}
||u_{0\eps}(\cdot)||_{L^p}$$ $$+ \int_0^t ||E_{n\eps}(t-\tau,
x-\cdot)||_{L^1} ||\nabla g_{\eps}(\theta u_{\eps})||_{L^{\infty}}
||u_{\eps}(\tau,\cdot)||_{L^p}d\tau.$$ Then, by (\ref{montreal})
$$
||u_{\eps}(t,\cdot)||_{L^p} \leq C|ln\eps|^{n(1-1/p)}+\int_0^t C
|ln\eps|^b ||u_{\eps}(\tau,\cdot)||_{L^p} d\tau.$$ Gronwall
inequality yields
$$
||u_{\eps}(t,\cdot)||_{L^p} \leq C |ln\eps|^{n(1-1/p)} \exp{(CT
|ln\eps|^b)}\leq C \eps^{-N}, $$ $\exists N>0,\; t\in [0,T), \;
T>0, \; x\in {\bf R}^n, \; b<1.$ The same holds for the first
derivative.

Consider $\beta^{th}$-derivative, $\beta \geq 2$,
$$
\partial_x^{\beta}u_{\eps}(t,x)=\int_{{\bf R}^n}
\partial_x^{\beta}E_{n\eps}(t,x-y) u_{0\eps}(y)dy+\int_0^t \int_{{\bf
R}^n}
\partial_x^{\beta}E_{n\eps}(t-\tau, x-y)
g_{\eps}(u_{\eps}(\tau,y))dyd\tau.$$ We have,
$$
||\partial_x^{\beta}u_{\eps}(t,\cdot)||_{L^p} \leq
||\partial_x^{\beta}E_{n\eps}(t, x-\cdot)||_{L^1}
||u_{0\eps}(\cdot)||_{L^p} $$ $$+\int_0^t ||\partial_x^{\beta}
E_{n\eps}(t-\tau, x-\cdot)||_{L^1} ||\nabla g_{\eps}(\theta
u_{\eps})||_{L^{\infty}} ||u_{\eps}(\tau,\cdot)||_{L^p}d\tau.$$
Then,
$$
||\partial_x^{\beta}u_{\eps}(t,\cdot)||_{L^p} \leq C
(ln|ln\eps|)^m |ln\eps|^{n(1-1/p)}+\int_0^t C (ln|ln\eps|)^m
|ln\eps|^{b} ||u_{\eps}(\tau, \cdot)||_{L^p}d\tau,$$ where
$m>\beta/2-1.$ Using the moderateness of $u_{\eps}(t,x)$ from the
first step we obtain
$$
||\partial_x^{\beta}u_{\eps}(t,\cdot)||_{L^p}\leq C
(ln|ln\eps|)^m |ln\eps|^{n(1-1/p)} +(CT (ln|ln\eps|)^m
 |ln\eps|^{b}\eps^{-N})\leq
C \eps^{-N},$$ $\exists N>0, \; t\in [0,T), \; x\in {\bf R}^n, \;
\eps<\eps_0,$ where $b<1,$ $m>\beta/2-1,$  $\beta\geq 2.$

Thus, $u_{\eps}(t,x)\in {\cal E}_{M,p}([0,T), {\bf R}^n).$

Consider the uniqueness. If we suppose that
$W_{\eps}(t,x)=u_{1\eps}(t,x)-u_{2\eps}(t,x)$ where
$u_{1\eps}(t,x)$ and $u_{2\eps}(t,x)$ are two solutions to the
equation (\ref{jedanaest}), then, we should to solve the equation
$$
\partial_tW_{\eps}(t,x)=\triangle W_{\eps}(t,x)+w_{\eps}(t,x)
W_{\eps}(t,x)+N_{\eps}(t,x)$$
$$
u_{0\eps}(x)=N_{0\eps}(x)\in {\cal N}_{p,q}({\bf R}^n), \;
N_{\eps}(t,x)\in {\cal N}_{C^1,(L^p, L^q)}([0,T), {\bf R}^n),$$
$$
w_{\eps}(t,x)=\int_0^t \nabla_ug_{\eps}(\theta
u_{1\eps}(t,x)+(1-\theta)u_{2\eps}(t,x))d\theta.$$ In integral
form we have
$$
W_{\eps}(t,x)=\int_{{\bf R}^n} E_{n\eps}(t,x-y)
N_{0\eps}(y)dy+\int_0^t \int_{{\bf R}^n} E_{n\eps}(t-\tau, x-y)
w_{\eps}(\tau,y) W_{\eps}(\tau,y)dyd\tau$$ $$+\int_0^t \int_{{\bf
R}^n} E_{n\eps}(t-\tau, x-y) N_{\eps}(\tau,y)dyd\tau.$$ In
$L^q$-norm, where $1\leq q \leq \infty$ we obtain
$$
||W_{\eps}(t,\cdot)||_{L^q} \leq ||E_{n\eps}(t,x-\cdot)||_{L^1}
||N_{0\eps}(\cdot)||_{L^q}+ \int_0^t ||E_{n\eps}(t-\tau,
x-\cdot)||_{L^1} ||w_{\eps}(\tau,\cdot)||_{L^{\infty}}$$
$$
||W_{\eps}(\tau,\cdot)||_{L^q}d\tau+\int_0^t ||E_{n\eps}(t-\tau,
x-\cdot) ||_{L^1} ||N_{\eps}(\tau,\cdot)||_{L^q}d\tau.$$ Then,
$$
||W_{\eps}(t,\cdot)||_{L^q} \leq  C  \eps^s+\int_0^t |ln\eps|^{b}
||W_{\eps}(\tau,\cdot)||_{L^q} d\tau+\int_0^t C  \eps^s d\tau.$$
Employ Gronwall inequality to obtain
$$
||W_{\eps}(t,\cdot)||_{L^q} \leq C \eps^s \exp{(CT
|ln\eps|^{b})}\leq C \eps^s,$$ $\forall s \in {\bf N}^n, \; x \in
{\bf R}^n, \; \eps<\eps_0,$ where $b<1.$ The same holds for the
first derivative.

For $\beta^{th}$-derivative, $\beta \geq 2,$
 we obtain
$$
\partial_x^{\beta}W_{\eps}(t,x)=\int_{{\bf R}^n}
\partial_x^{\beta}E_{n\eps}(t,x-y) N_{0\eps}(y)dy+\int_0^t
\int_{{\bf R}^n} \partial_x^{\beta}E_{n\eps}(t-\tau,x-y)
w_{\eps}(\tau,y)$$ $$W_{\eps}(\tau,y) dyd\tau+ \int_0^t \int_{{\bf
R}^n} \partial_x^{\beta}E_{n\eps}(t-\tau,x-y) N_{\eps}(\tau,y)
dyd\tau.$$ In $L^q$-norm we have
$$
||\partial_x^{\beta}W_{\eps}(t,\cdot)||_{L^q} \leq C
||\partial_x^{\beta}E_{n\eps}(t,x-\cdot)||_{L^1}
||N_{0\eps}(\cdot)||_{L^q}+ \int_0^t
||\partial_x^{\beta}E_{n\eps}(t-\tau, x-\cdot)||_{L^1}$$
$$
||w_{\eps}(\tau,\cdot)||_{L^{\infty}}
||W_{\eps}(\tau,\cdot)||_{L^q}d\tau+ \int_0^t
||\partial_x^{\beta}E_{n\eps}(t-\tau, x-\cdot)||_{L^1}
||N_{\eps}(\tau,\cdot)||_{L^q} d\tau,$$ for $m>\beta/2-1.$ By the
first step and (\ref{montreal})
$$
||\partial_x^{\beta} W_{\eps}(t,\cdot)||_{L^q} \leq
C(ln|ln\eps|)^m \eps^s +\int_0^t C (ln|ln\eps|)^{m} |ln\eps|^b
\eps^s d\tau+\int_0^t C (ln|ln\eps|)^m \eps^s d\tau$$ $ \leq C
\eps^s,\;\forall s \in {\bf N}, \; t\in [0,T), \; x\in {\bf R}^n,
\; \eps<\eps_0,$ under the condition $b<1$,  where $m>\beta/2-1.$

Finally, $W_{\eps}(t,x)\in {\cal N}_{C^1,(L^p,L^q)}([0,T), {\bf
R}^n)$, $1\leq p,q \leq \infty.$ Thus, the solution is unique in
the Colombeau vector space ${\cal G}_{C^1,(L^p,L^q)}([0,T), {\bf
R}^n)$ for those choices $1\leq p,q\leq \infty$ for which
corresponding Colombeau vector space is an algebra with
multiplication. This is the consequence of the Sobolev imbedding
theorems, cf. \cite{adams}.

$\Box$

\subsection{\bf Parabolic equation with nonlinear conservative
term (\ref{studenti})}

The equation (\ref{studenti}) has the following regularized
integral form
\begin{equation}\label{markic}
u_{\eps}(t,x)=E_{n\eps}(t,x)*\mu_{\eps}(x)+\int_0^t \nabla
E_{n\eps}(t-\tau,
\cdot)*\vec{g}_{\eps}(u_{\eps}(\tau,\cdot))d\tau,\;
\mu_{\eps}(x)=\delta_{\eps}(x),
\end{equation}
where the regularization for each term with subscript $\eps$
except for the heat kernel,
 are
given in (\ref{dvanaest}).

We shall regularize the gradient of the heat kernel as follows: $
\nabla E_{n\eps}(t,x)=$ $\nabla E_n (t,x)$ $*$ $\phi_{\eps}(t). $
In $L^1$-norm, we have
$$
||\nabla E_{n\eps}(t-\tau, x-\cdot)||_{L^1}=\int_{{\bf R}^n}
|\int_{\bf R} \nabla E_n(t-\tau,x-y) \phi_{\eps}(\tau)d\tau|dy$$
$$\leq C \int_{{\bf R}^n} |\int_{\bf R} (t-\tau)^{-n/2-1} \exp{(-
|x-y|^2)/(4(t-\tau))}\phi_{\eps}(\tau)d\tau| dy.$$ Then,
$$
||\partial_x^{\beta}\nabla E_{n\eps}(t-\tau, x-\cdot)||_{L^1}=
\int_{{\bf R}^n}|\int_{\bf R}
\partial_y^{\beta}\nabla E_{n}(t-\tau, x-y)
\phi_{\eps}(\tau)d\tau| dy $$
$$\leq
C \int_{{\bf R}^n} \int_{\bf R} |(t-\tau)^{-n/2-1-\beta}|
|\exp{(-|x-y|^2/(4(t-\tau)))}| |\phi_{\eps}(\tau)|d\tau dy.$$
Apply Fubini theorem to obtain
$$
||\partial_x^{\beta}\nabla E_{n\eps}(t-\tau, x-\cdot)||_{L^1}\leq
C \int_{\bf R}\int_{{\bf R}^n} |(t-\tau)^{-n/2-\beta-1}|
|\phi_{\eps}(\tau)|$$ $$  |\exp{(-|x-y|^2/(4(t-\tau)))}|dy
d\tau.$$ Then, we continue
$$
\leq C \int_{\bf R} |(t-\tau)^{-\beta/2-1/2}
||(t-\tau)^{(\beta+1)/2}\partial_x^{\beta}\nabla E_n(t-\tau,
x-\cdot)||_{L^1} d\tau$$ $$\leq C \int_{\bf R}
(t-\tau)^{-\beta/2-1/2} \phi_{\eps}(\tau) d\tau.$$ Apply Lemma
\ref{lemmaone} to obtain
\begin{equation}\label{evorucak}
||\partial_x^{\beta}\nabla E_{n\eps}(t-\tau, x-\cdot)||_{L^1}\leq
\left\{\begin{array}{ll} C \; \; \mbox{for} \; \beta <1\\
C(ln|ln\eps|)^m \; \mbox{for} \; m>(\beta-1)/2, \; \beta \geq
1.\end{array}\right.\end{equation}

The initial data are regularized by delta sequence, nonlinear term
$g(u)\in L_{loc}^{\infty}([0,T),$ $ {\bf R}^n)$ which is nonlinear
and non-Lipschitz, it
 is regularized by cut-off and
the fractional part of the heat kernel is regularized as described
in Lemma \ref{lemmaone}. We have the following Theorem.
\begin{th}
 The regularized
equation (\ref{markic}) have a unique solution $[u_{\eps}] \in
{\cal G}_{C^1, (L^p,L^q)}$ $([0,T),$ $ {\bf R}^n)$ for those
choices $1\leq p,q \leq \infty$ for which the corresponding
Colombeau space is an algebra with multiplication. In a limiting
case, when $\eps \to 0$, we obtain the solution to the equation
(\ref{studenti}).
\end{th}
{\bf Proof.}  Taking the $L^p$-norm, $1\leq p \leq \infty$ of
(\ref{markic}) we obtain
$$
||u_{\eps}(t,\cdot)||_{L^p} \leq ||E_{n\eps}(t,\cdot)||_{L^1}
||\mu_{\eps}(\cdot)||_{L^p}$$ $$+\int_0^t ||\nabla
E_{n\eps}(t-\tau, x-\cdot)||_{L^1} ||\nabla g_{\eps}(\theta
u_{\eps})||_{L^{\infty}} ||u_{\eps}(\tau, \cdot)||_{L^p} d\tau.$$
Since $\mu_{\eps}(x)=\delta_{\eps}(x)=|ln\eps|^{an} \phi(x \cdot
|ln\eps|)$ and by $(|ln\eps|)$-boundedeness of gradient we have
$$
||u_{\eps}(t,\cdot)||_{L^p} \leq C  |ln\eps|^{n(a-1/p)}+\int_0^t C
|ln\eps|^{b} ||u_{\eps}(\tau,\cdot)||_{L^p} d\tau.$$  By Gronwall
inequality
$$
||u_{\eps}(t,\cdot)||_{L^p} \leq C |ln\eps|^{n(a-1/p)} \exp{(CT
|ln\eps|^{b})}\leq C \eps^{-N},$$ $\exists N>0, \; t\in [0,T), \;
x\in {\bf R}^n, \; b<1, \; a>0.$
\begin{rem}
If we use $(ln|ln\eps|)$-boundedeness for $|\nabla g_{\eps}|$ and
for the heat semigroup then, there is no condition on parameters.
Slower growth of mollifiers contributes to the moderateness.
\end{rem}
Suppose that $\beta \in {\bf N}_0^n,$ $\beta >0.$
 Then, we have from
(\ref{markic})
$$\partial_x^{\beta}u_{\eps}(t,x)=\partial_x^{\beta}E_{n\eps}(t,x)*\mu_{\eps}(x)+\int_0^t
\int_{{\bf R}^n}
\partial_x^{\beta}\nabla E_{n\eps}(t-\tau, x-y)
\vec{g}_{\eps}(u_{\eps}(\tau,y))dy d\tau.$$  By (\ref{evorucak})
$$
||\partial_x^{\beta} u_{\eps}(t,\cdot)||_{L^1} \leq C
(ln|ln\eps|)^m |ln\eps|^{n(1-1/p)} + \int_0^t C (ln|ln\eps|)^{m}
||\nabla g_{\eps}(\theta u_{\eps})||_{L^{\infty}}
||u_{\eps}(\tau,\cdot)||_{L^p} d\tau.$$ Using the moderateness of
$u_{\eps}(t,x)$ to obtain
$$
||\partial_x^{\beta}u_{\eps}(t,\cdot)||_{L^p} \leq C
(ln|ln\eps|)^m |ln\eps|^{n(1-1/p)} + CT (ln|ln\eps|)^m
|ln\eps|^{b} \eps^{-N}\leq C \eps^{-N},$$ $ \exists N>0, \; t\in
[0,T), \; x\in {\bf R}^n, \; \eps<\eps_0, \; m>(\beta-1)/2$ under
the condition $b<1$ what is included in the condition imposed at
the first step.

Thus, we have the moderateness in the space $u_{\eps}(t,x) \in
{\cal E}_{M,p}([0,T), {\bf R}^n).$

Concerning the uniqueness we should to  solve the equation
$$\partial_t W_{\eps}(t,x)=\triangle W_{\eps}(t,x)+w_{\eps}(t,x)
W_{\eps}(t,x)+ N_{\eps}(t,x)$$
$$W_{0\eps}(0,x)=N_{0\eps}(x) \in {\cal N}_{ p,q}({\bf
R}^n)$$ where $1\leq p,q \leq \infty,$ $ N_{\eps}(t,x)\in  {\cal
N}_{C^1, (L^p,L^q)}([0,T), {\bf R}^n)$,
$W_{\eps}(t,x)=u_{1\eps}(t,x)-u_{2\eps}(t,x)$ where
$u_{1\eps}(t,x)$ and $u_{2\eps}(t,x)$ are two solutions to the
regularized  equation to (\ref{studenti}) and
$w_{\eps}(t,x)=\int_0^t \nabla g_{\eps}(t, \theta
u_{1\eps}(t,x)+(1-\theta)u_{2\eps})d\theta,$ $0<\theta <1.$
 In integral
form we have
$$
W_{\eps}(t,x)=E_{n\eps}(t,x)*N_{0\eps}(x)+\int_0^t \int_{{\bf
R}^n} \nabla E_{n\eps}(t-\tau,x-y) w_{\eps}(\tau,y)
W_{\eps}(\tau,y)dy d\tau$$ $$+\int_0^t \int_{{\bf R}^n} \nabla
E_{n\eps}(t-\tau, x-y) N_{\eps}(\tau,y) dy d\tau.$$ In $L^q$-norm
where $1\leq q \leq \infty$, we have
$$
||W_{\eps}(t,\cdot)||_{L^q} \leq ||E_{n\eps}(t,\cdot)||_{L^1}
||N_{0\eps}(\cdot)||_{L^q}+\int_0^t ||\nabla E_{n\eps}(t-\tau,
x-\cdot)||_{L^1} ||w_{\eps}(\tau,\cdot)||_{L^{\infty}}$$
$$
||W_{\eps}(\tau,\cdot)||_{L^q} d\tau+ \int_0^t ||\nabla
E_{n\eps}(t-\tau, x-\cdot)||_{L^1} ||N_{\eps}(\tau,\cdot)||_{L^q}
d\tau.$$ Then, we obtain
$$
||W_{\eps}(t,\cdot)||_{L^q} \leq C  \eps^s+\int_0^t C |ln\eps|^{b}
||W_{\eps}(\tau,\cdot)||_{L^q} d\tau.$$ Employ Gronwall inequality
to obtain
$$
||W_{\eps}(t,\cdot)||_{L^q} \leq C  \eps^s \exp{(CT
|ln\eps|^{b})}\leq C \eps^s,$$ $ \forall s \in {\bf N}, \; t\in
[0,T), \; x\in {\bf R}^n, \; \eps<\eps_0,$ where $b<1$.

For $\beta^{th}$-derivative we have for $\beta >0$
$$
\partial_x^{\beta}W_{\eps}(t,x)=\partial_x^{\beta}E_{n\eps}(t,x)*N_{0\eps}(x)+
\int_0^t \int_{{\bf R}^n}
\partial_x^{\beta} \nabla E_{n\eps}(t-\tau, x-y)
w_{\eps}(\tau,y) W_{\eps}(\tau,y) dy d\tau$$ $$+\int_0^t
\int_{{\bf R}^n} \partial_x^{\beta} \nabla E_{n\eps}(t-\tau, x-y)
N_{\eps}(\tau,y) dy d\tau.$$ In $L^q$-norm we have
$$
||\partial_x^{\beta}W_{\eps}(t,\cdot)||_{L^q} \leq
||\partial_x^{\beta}E_{n\eps}(t,\cdot)||_{L^1}
||N_{0\eps}(\cdot)||_{L^q}+\int_0^t ||\partial_x^{\beta} \nabla
E_{n\eps}(t-\tau, x-\cdot)||_{L^1}$$
$$
||w_{\eps}(\tau,\cdot)||_{L^{\infty}}
||W_{\eps}(\tau,\cdot)||_{L^q}d\tau+\int_0^t ||\partial_x^{\beta}
\nabla E_{n\eps}(t-\tau, x-\cdot)||_{L^1}
||N_{\eps}(\tau,\cdot)||_{L^q}d\tau,$$ and for $m>(\beta-1)/2$ due
to null property of $W_{\eps}(t,x)$ obtained at the first step of
induction we obtain
$$
||\partial_x^{\beta}W_{\eps}(\tau,\cdot)||_{L^q} \leq C
(ln|ln\eps|)^{m} \eps^s+\int_0^t C (ln|ln\eps|)^m |ln\eps|^{b}
\eps^s d\tau+\int_0^t C (ln|ln\eps|)^{m}\eps^s d\tau $$ $\leq C
\eps^s, \; \forall s \in {\bf N}, \; t\in [0,T), \; x\in {\bf
R}^n, \; \eps <\eps_0,$ $m>(\beta-1)/2$,
 under the condition imposed at the first
step, $b<1.$

Thus, $W_{\eps}(t,x) \in {\cal N}_{C^1, (L^p,L^q)}([0,T), {\bf
R}^n),$ $1\leq p,q \leq \infty.$ Follows, there exists  the unique
solution to the equation (\ref{studenti}) in the space ${\cal
G}_{C^1, (L^p, L^q)}([0,T), {\bf R}^n),$ $1\leq p,q \leq \infty$,
for those  choices of $p,q$ for which corresponding  Colombeau
space is an algebra with multiplication. If $p=q=2$ due to the
Sobolev imbedding theorem for $n \leq 3,$ $L^2({\bf R}^n) \subset
L^{\infty}({\bf R}^n)$ and the space ${\cal
G}_{C^1,(L^2,L^2)}([0,T), {\bf R}^n)$ is an algebra with
multiplication (cf. \cite{n+p+r}). $\Box$

\subsection{\bf Nonlinear parabolic equation with Schr\"odinger
kernel (\ref{zadaci})}

Consider nonlinear parabolic equation with Schr\"odinger kernel
(\ref{zadaci}) where $u_0(x)=\delta(x),$ $V(x)=\delta(x),$
$g(u)\in L_{loc}^{\infty}([0,T), {\bf R}^n)$ is nonlinear,
non-Lipschitz and it is regularized as in \cite{polar} to obtain
$$|\nabla g_{\eps}(\theta u_{\eps})|\leq C |ln\eps|^b, \;
(\mbox{resp.} \; (ln|ln\eps|)).$$ Initial data and potential are
regularized as follows:
$$
u_{0\eps}(x)=|ln\eps|^{an} \phi(x \cdot |ln\eps|), \;
V_{\eps}(x)=|ln\eps|^{cn} \phi(x \cdot |ln\eps|),$$ and the heat
semigroup has the regularization   described in Section
\ref{petjedan}.

We have in integral form
\begin{equation}\label{markovici}
u_{\eps}(t,x)=E_{n\eps}(t,x)*\mu_{\eps}(x)+\int_0^t \int_{{\bf
R}^n} E_{n\eps}(t-\tau, x-y) V_{\eps}(y) u_{\eps}(\tau,y) dy
d\tau\end{equation}
$$+ \int_0^t \int_{{\bf R}^n} E_{n\eps}(t-\tau, x-y)
g_{\eps}(u_{\eps}(\tau,y)) dy d\tau.$$
 Then, in
$L^p$-norm, $1\leq p \leq \infty,$ we have
$$
||u_{\eps}(t,\cdot)||_{L^p} \leq ||E_{n\eps}(t,\cdot)||_{L^1}
||\mu_{\eps}(\cdot)||_{L^p}+\int_0^t
||E_{n\eps}(t-\tau,x-\cdot)||_{L^1}
||V_{\eps}(\cdot)||_{L^{\infty}} ||u_{\eps}(\tau,\cdot)||_{L^p}
d\tau$$ $$+\int_0^t ||E_{n\eps}(t-\tau,x-\cdot)||_{L^1} ||\nabla
g_{\eps}(\theta u_{\eps})||_{L^{\infty}}
||u_{\eps}(\tau,\cdot)||_{L^p} d\tau.$$ By regularization
(\ref{montreal}) given for heat kernel
 we have
for $\beta<1$
$$
||u_{\eps}(t,\cdot)||_{L^p} \leq C |ln\eps|^{n(1-1/p)}+\int_0^t C
|ln\eps|^{cn} ||u_{\eps}(\tau,\cdot)||_{L^p} d\tau$$ $$+\int_0^t C
|ln\eps|^{b} ||u_{\eps}(\tau,\cdot)||_{L^p} d\tau.$$ Employ
Gronwall inequality to obtain
$$
||u_{\eps}(t,\cdot)||_{L^p} \leq C |ln\eps|^{n(1-1/p)} \exp{(CT
(|ln\eps|^{cn}+|ln\eps|^{b}))} \leq C \eps^{-N},$$ $\exists N>0,
\; \eps< \eps_0, \; t\in [0,T), \; x\in {\bf  R}^n, \; cn<1, \;
b<1.$ The same holds for the first derivative.
\begin{rem}
The condition
$\left\{\begin{array}{ll} cn<1\\
b<1
\end{array}\right.$
gives the link
 between the singularity  of the heat semigroup, nonlinearity of
 $g(u)$ and the singularity of the potential.
 Singularity of the initial data is given by
  $n(1-1/p)>0$ what  holds and it is necessary that this part is
 logarithmically bounded.
\end{rem}
Consider $\beta^{th}$-derivative, $\beta \geq 2$. We have in
(\ref{markovici})
$$
\partial_x^{\beta}u_{\eps}(t,x)=
\partial_x^{\beta}E_{n\eps}(t,x)*\mu_{\eps}(x)+\int_0^t \int_{{\bf
R}^n}
\partial_x^{\beta}E_{n\eps}(t-\tau, x-y) V_{\eps}(y)
u_{\eps}(\tau,y)dy d\tau$$ $$+\int_0^t \int_{{\bf R}^n}
\partial_x^{\beta} E_{n\eps}(t-\tau, x-y)
g_{\eps}(u_{\eps}(\tau,y))dy d\tau.$$ In $L^p$-norm, where $1\leq
p \leq \infty$, we obtain
$$
||\partial_x^{\beta}u_{\eps}(t,\cdot)||_{L^p}\leq
||\partial_x^{\beta} E_{n\eps}(t,\cdot)||_{L^1}
||\mu_{\eps}(\cdot)||_{L^p}+\int_0^t ||\partial_x^{\beta}
E_{n\eps}(t-\tau, x-\cdot)||_{L^1}
||V_{\eps}(\cdot)||_{L^{\infty}}$$
$$
||u_{\eps}(\tau,\cdot)||_{L^p}d\tau+ \int_0^t ||\partial_x^{\beta}
E_{n\eps}(t-\tau, x-\cdot)||_{L^1} ||\nabla g_{\eps}(\theta
u_{\eps})||_{L^{\infty}} ||u_{\eps}(\tau,\cdot)||_{L^p} d\tau,$$
and by the regularization we have for $m>\beta/2-1$
$$
||\partial_x^{\beta}u_{\eps}(t,\cdot)||_{L^p} \leq C
(ln|ln\eps|)^m |ln\eps|^{n(1-1/p)}+\int_0^t C (ln|ln\eps|)^b
|ln\eps|^{cn} \eps^{-N} d\tau$$ $$+\int_0^t C (ln|ln\eps|)^m
|ln\eps|^{b} \eps^{-N}d\tau.$$ We used the moderateness of
$u_{\eps}(t,x)$ obtained at the first step. Thus,
$$||\partial_x^{\beta}(t,\cdot)||_{L^p}\leq C \eps^{-N}, \;
\exists N>0, \; t\in [0,T), \; x\in {\bf R}^n, \; \eps<\eps_0,$$
under the conditions $b<1,$ $cn<1,$ imposed at the first step, and
$m>\beta/2-1.$

Concerning the uniqueness we should to solve the equation
$$
\partial_tW_{\eps}(t,x)=\triangle W_{\eps}(t,x)+V_{\eps}(x)
W_{\eps}(t,x)+w_{\eps}(t,x) W_{\eps}(t,x)+N_{\eps}(t,x)$$
$$ W_{\eps}(0,x)=N_{0\eps}(x)\in {\cal N}_{ p, q}({\bf
R}^n),$$ where $N_{\eps}(t,x) \in {\cal N}_{C^1, (L^p,L^q)}([0,T),
{\bf R}^n),$ $W_{\eps}(t,x)=u_{1\eps}(t,x)-u_{2\eps}(t,x)$ where
$u_{1\eps}(t,x)$ and $u_{2\eps}(t,x)$
 are two solutions to the regularization for the equation (\ref{zadaci}),
 $w_{\eps}(t,x)=\int_0^1 \nabla g_{\eps}(t, \theta
 u_{1\eps}(t,x)+(1-\theta) u_{2\eps}(t,x))d\theta.$

 In integral form we have
\begin{equation}\label{milanka}
W_{\eps}(t,x)=E_{n\eps}(t,x)*N_{0\eps}(x)+\int_0^t \int_{{\bf
R}^n} E_{n\eps}(t-\tau, x-y) V_{\eps}(y) W_{\eps}(\tau,y) dy
d\tau\end{equation} $$+ \int_0^t \int_{{\bf R}^n}
E_{n\eps}(t-\tau, x-y) g_{\eps}(u_{\eps}(\tau,y))dy d\tau
+\int_0^t \int_{{\bf R}^n} E_{n\eps}(t-\tau, x-y)
N_{\eps}(\tau,y)dy d\tau.
$$
In $L^q$-norm, $1\leq q \leq \infty,$ we obtain
$$
||W_{\eps}(t,\cdot)||_{L^q} \leq C ||E_{n\eps}(t,\cdot)||_{L^1}
||N_{0\eps}(\cdot)||_{L^q}+\int_0^t ||E_{n\eps}(t-\tau,
x-\cdot)||_{L^1} ||V_{\eps}(\cdot)||_{L^{\infty}}$$
$$
||W_{\eps}(\tau,\cdot)||_{L^q} d\tau + \int_0^t
||E_{n\eps}(t-\tau, x-\cdot)||_{L^1}
||w_{\eps}(\tau,\cdot)||_{L^{\infty}}
||W_{\eps}(\tau,\cdot)||_{L^q}d\tau $$ $$+\int_0^t
||E_{n\eps}(t-\tau, x-\cdot)||_{L^1}
||N_{\eps}(\tau,\cdot)||_{L^q}d\tau.$$ Employ the regularizations
to obtain
$$
||W_{\eps}(t,\cdot)||_{L^q} \leq C  \eps^s+\int_0^t C
|ln\eps|^{cn} ||W_{\eps}(\tau,\cdot)||_{L^q} d\tau+\int_0^t C
|ln\eps|^{b} ||W_{\eps}(\tau,\cdot)||_{L^q} d\tau.$$  By Gronwall
inequality
$$
||W_{\eps}(t,\cdot)||_{L^q} \leq C  \eps^s
\exp{(CT(|ln\eps|^{cn}+|ln\eps|^{b})})\leq C \eps^s,$$ $\forall s
\in {\bf N}, \; x\in {\bf R}^n, \; t\in [0,T), \; \eps<\eps_0,$
where the condition $c n <1,$ $b<1$ is imposed at the first step
in the proof of moderateness. The same holds for the first
derivative due to (\ref{montreal}).

For $\beta^{th}$-derivative, $ \beta \geq 2,$ we have in
(\ref{milanka})
$$
\partial_x^{\beta}W_{\eps}(t,x)= \partial_{x}^{\beta}
E_{n\eps}(t,x)* N_{0\eps}(x)+\int_0^t \int_{{\bf R}^n}
\partial_x^{\beta} E_{n\eps}(t-\tau, x-y)
V_{\eps}(y) W_{\eps}(\tau,y) dy d\tau$$ $$+\int_0^t \int_{{\bf
R}^n}
\partial_x^{\beta} E_{n\eps}(t-\tau,x-y) w_{\eps}(\tau,y)
W_{\eps}(\tau,y) dy d\tau +\int_0^t \int_{{\bf R}^n}
\partial_x^{\beta} E_{n\eps}(\tau,y) N_{\eps}(\tau,y)dy d\tau.$$
 Then, in $L^q$-norm, $1\leq q \leq \infty$, we have
$$
||\partial_x^{\beta}W_{\eps}(t,\cdot)||_{L^q} \leq
||\partial_x^{\beta} E_{n\eps}(t,\cdot)||_{L^1}
||N_{0\eps}(\cdot)||_{L^q}+\int_0^t
||\partial_x^{\beta}E_{n\eps}(t-\tau, x-\cdot)||_{L^1}
||V_{\eps}(\cdot)||_{L^{\infty}}$$ $$
||W_{\eps}(\tau,\cdot)||_{L^q} d\tau+\int_0^t ||\partial_x^{\beta}
E_{n\eps}(t-\tau, x-\cdot)||_{L^1}
||w_{\eps}(\tau,\cdot)||_{L^{\infty}}
||W_{\eps}(\tau,\cdot)||_{L^q} d\tau $$ $$+\int_0^t
||\partial_x^{\beta} E_{n\eps}(\tau, \cdot)||_{L^1}
||N_{\eps}(\tau,\cdot)||_{L^q} d\tau.$$ Setting the regularization
(\ref{montreal}) and due to the null property of $W_{\eps}(t,x)$
we obtain
$$
||\partial_x^{\beta}W_{\eps}(t,\cdot)||_{L^q} \leq C
(ln|ln\eps|)^{m} \eps^s+\int_0^t C (ln|ln\eps|)^m |ln\eps|^{cn}
\eps^{s} d\tau$$ $$+\int_0^t C (ln|ln\eps|)^m |ln\eps|^{b}
\eps^{s} d\tau\leq C \eps^{s},$$ $ \forall s \in {\bf N}, \; \eps<
\eps_0, \; t\in [0,T), \; x\in {\bf R}^n,$ $m>\beta/2-1,$
 under the condition
$b<1$, $cn<1$, imposed at the first step of induction.

Thus,$W_{\eps}(t,x) \in {\cal N}_{C^1, (L^p, L^q)}([0,T), {\bf
R}^n)$, $ 1\leq p,q \leq \infty$.

Follows, the solution is unique in the spaces ${\cal
G}_{C^1,(L^p,L^q)}([0,T), {\bf R}^n)$, $1\leq p,q \leq \infty,$
for those choices of  $p,q$ for which the corresponding Colombeau
space is an algebra with multiplication. That is the consequence
of Sobolev imbedding theorems. $\Box$
\begin{rem}
We can remain the heat semigroup to stay $|ln\eps|$-bounded and
set the potential and $|\nabla g_{\eps}|$ to be
$(ln|ln\eps|)$-bounded to obtain existence-uniqueness theorems in
corresponding Colombeau algebras  with appropriate conditions on
parameters what leads to the balance between the singularities of
the initial data, heat semigroup, potential and the nonlinearity
of $g.$
\end{rem}

\section{\bf Application to linear Schr\"odinger equation with
singular potential and initial data}

Consider the linear Schr\"odinger equation with strongly singular
potential and initial data
\begin{equation}\label{donau}
\frac{1}{i}
 \partial_tu(t,x)=(\triangle -V(x))u(t,x), \;
u(0,x)=u_0(x)=\delta(x), \; V(x)=\delta(x).
\end{equation}
The following regularization for delta distribution will be used:
$$
u_{0\eps}(x)=|ln\eps|^{an}\phi(x \cdot |ln\eps|), \;
V_{\eps}(x)=|ln\eps|^{cn} \phi(x\cdot |ln\eps|), \; (\mbox{resp.
using } \; (ln|ln\eps|)),$$ where $\phi(x) \in C_0^{\infty}({\bf
R}^n),$
 $\phi(x)\geq 0,$ $\int
\phi(x)dx=1.$

 We shall give the regularization for
Schr\"odinger  semigroup  with respect to $t$ by delta sequence to
handle the existence-uniqueness result in the vector space ${\cal
G}_{C^1, (L^p, L^q)}$ $([0,T),$ $ {\bf R}^n).$ Denote by
$S_{n\eps}(t,x)=S_n(t,x)*\phi_{\eps}(t),$ where $S_n(t,x)=(4\pi
t)^{-n/2},$ $  \exp{(i |x|^2/(4t)}.$
 Then,
$$
\sup_x |S_{n\eps}(t,x)| \leq  \sup_{x}|\int_0^t S_n(t-\tau,y)
\phi_{\eps}(\tau)d\tau|$$ $$\leq \sup_x \int_0^t |(4\pi
(t-\tau))^{-n/2}| |\exp{(i|x|^2/(4(t-\tau))})| |\phi_{\eps}(\tau)|
d\tau \leq C \int_0^t |(t-\tau)^{-n/2}| |\phi_{\eps}(\tau)|
d\tau,$$ since $\sup_x |\exp{(i|x|^2/(4(t-\tau)))} | \leq 1.$ This
is the fractional derivative of $\delta$-sequence and by Lemma
\ref{lemmaone} it follows
\begin{equation}\label{mario}
\sup_{x} |S_{n\eps}(t,x)| \leq \left\{\begin{array}{ll}C \; \; \;
\mbox{for} \; n<2\\
C(ln|ln\eps|)^m, \; \mbox{for} \; m> n/2-1, \; n\geq 2.
\end{array}\right.
\end{equation}
Consider the $\beta^{th}$-derivative, $\beta \in {\bf N}_0^n,$
 of Schr\"odinger semigroup. We have $
\partial_x^{\beta}S_{n\eps}(t,x)=\partial_x^{\beta}S_n(t,x)*\phi_{\eps}(t).$
Since $\sup_x |\partial_x^{\beta} S_{n\eps}(t,x)| \leq
C(t-\tau)^{-n/2-\beta}$ we obtain
$$
\sup_x |\partial_x^{\beta} S_{n\eps}(t,x)| \leq  C \int_0^t
(t-\tau)^{-n/2-\beta} \phi_{\eps}(\tau) d\tau.$$ By Lemma
\ref{lemmaone}
\begin{equation}\label{kale}
\sup_x |\partial_x^{\beta} S_{n\eps}(t,x)| \leq
\left\{\begin{array}{ll} C \; \; \mbox{for} \; n=1, \; \beta=0\\
C(ln|ln\eps|)^m \; \mbox{for} \; m>\beta+n/2-1, \; \beta \geq 0,
\; n\geq 2.
\end{array}\right.\end{equation}
We indicate interpolation spaces for noninteger $m$, from
\cite{temam}. Suppose that X,Y are two Hilbert space, $X \subset
Y,$ $X$ is dense in $Y$ and the injection is continuous. The
interpolation between the spaces provides a family of Hilbert
spaces $[X,Y]_{\theta},$ $0\leq \theta \leq 1,$ with the following
properties: $[X,Y]_0=X,$ $[X,Y]_1=Y$ and $X \subset [X,Y]_{\theta}
\subset Y$, where the injection is continuous and each space is
dense in succeeding one. Then,
$$
||f||_{[X,Y]_{\theta}} \leq c(\theta) ||f||_{X}^{1-\theta}
||f||_{Y}^{\theta}, \; \forall f \in X, \; \forall \theta, \;
1\leq\theta \leq 1.$$ We apply this to obtain norm in
$L^{p/(p+1)}$-space. Since $\frac{p}{p+1}=1-\frac{1}{p+1}$ then
this space is settled
 between $L^2$ and $L^1.$ We have $L^2 \subset [L^2,L^1]_{\theta}
\subset L^1,$ $1\leq \theta \leq 1.$ Then,
$$
||u_{0\eps}(\cdot)||_{L^{p/(p+1)}} \leq c(\theta)
||u_{0\eps}(\cdot)||_{L^1}^{1-\theta}
||u_{0\eps}(\cdot)||_{L^2}^{\theta}, \; 0 \leq \theta \leq 1$$ and
$$
||u_{0\eps}(\cdot)||_{L^{p/(p+1)}}\leq c(\theta)
|ln\eps|^{n(a-1)(1-\theta)+n(a-1/2)\theta}, \; 1\leq \theta \leq
1.$$ Thus,
\begin{equation} \label{mleko}
||u_{0\eps}(\cdot)||_{L^{p/(p+1))}} \leq c(\theta)
|ln\eps|^{n(a-1+\theta/2)}, \; \mbox{where} \; 0\leq \theta \leq
1.
\end{equation}
Now we can prove the following theorem.
\begin{th}
Regularized equation to Schr\"odinger equation (\ref{donau}) with
strongly singular initial data and potential where the
regularizations for potential and initial data are given by
(\ref{dvanaest}) and the regularization for Schr\"odinger
semigroup is described above, has a unique solution in the space
${\cal G}_{C^1, (L^p,L^q)}([0,T), {\bf R}^n)$, $1\leq p,q \leq
\infty,$ for those choices $p,q$ for which corresponding Colombeau
vector space is an algebra with multiplication (for example, for
$p=q=2,$ $n\leq 3,$
 this holds).
\end{th}
{\bf Proof.} The  regularized  equation for (\ref{donau}) has the
following integral form
\begin{equation}\label{akademiceskaja}
u_{\eps}(t,x)=S_{n\eps}(t,x)*u_{0\eps}(x)+\int_0^t \int_{{\bf
R}^n} S_{n\eps}(t-\tau,x-y) V_{\eps}(y) u_{\eps}(\tau,y)dy d\tau.
\end{equation}
In $L^p$-norm we have
$$
||u_{\eps}(t,\cdot)||_{L^p} \leq
||S_{n\eps}(t,\cdot)||_{L^{\infty}}
||u_{0\eps}(\cdot)||_{L^{p/(p+1)}}$$ $$ +\int_0^t \int_{{\bf R}^n}
||S_{n\eps}(t-\tau, x-\cdot)||_{L^{\infty}}
||V_{\eps}(\cdot)||_{L^1} ||u_{\eps}(\tau,\cdot)||_{L^p} d\tau.$$
Since $||V_{\eps}(\cdot)||_{L^1} \leq C |ln\eps|^{n(c-1)}$ and by
(\ref{mleko}) and (\ref{mario}) we obtain
$$
||u_{\eps}(t,\cdot)||_{L^p}\leq   \left\{ \begin{array}{ll} C
\; \; n<2\\
C(ln|ln\eps|)^m, \; n \geq 2, \;  m>n/2-1\end{array}\right.
|ln\eps|^{n(a-1+\theta/2)}$$
$$+
\int_0^t
 \left\{\begin{array}{ll}
C \; \; \; n<2\\
C(ln|ln\eps|)^m, \; n \geq 2, \; m>n/2-1 \end{array}\right.
|ln\eps|^{n(c-1)} ||u_{\eps}(\tau,\cdot)||_{L^p} d\tau.$$ By
Gronwall inequality
$$
||u_{\eps}(t,\cdot)||_{L^p} \leq  \left\{ \begin{array}{ll} C
\; \; n<2\\
C(ln|ln\eps|)^m, \; n \geq 2, \;  m>n/2-1\end{array}\right.
|ln\eps|^{n(a-1+\theta/2)}
$$
$$
\exp{(CT \left\{\begin{array}{ll}
C \; \; \; n<2\\
C(ln|ln\eps|)^m, \; n \geq 2, \; m>n/2-1 \end{array}\right.
|ln\eps|^{n(c-1)})}.$$ Thus,
$$
||u_{\eps}(t,\cdot)||_{L^p} \leq C \eps^{-N}, \; \exists N>0, \;
x\in {\bf R}^n, \; t\in [0,T), \; \eps<\eps_0, \; 1\leq p \leq
\infty, \; c<1+1/n.$$ Consider $\beta^{th}$-derivative, $\beta \in
{\bf N}_0^n,$ $\beta \geq 1,$
 of (\ref{akademiceskaja})
$$
\partial_x^{\beta}u_{\eps}(t,x)=
\partial_x^{\beta} S_{n\eps}(t,x)*u_{0\eps}(x)+\int_0^t \int_{{\bf
R}^n}
\partial_x^{\beta} E_{n\eps}(t-\tau, x-y) V_{\eps}(y)
u_{\eps}(\tau,y)dy d\tau.$$ Then, we have
$$
||\partial_x^{\beta}(t,\cdot)||_{L^p} \leq
||\partial_x^{\beta}S_{n\eps}(t,\cdot)||_{L^{\infty}}
||u_{0\eps}(\cdot)||_{L^{p/(p+1)}}$$ $$+\int_0^t \int_{{\bf R}^n}
||\partial_x^{\beta}S_{n\eps}(t-\tau, x-\cdot)||_{L^{\infty}}
||V_{\eps}(\cdot)||_{L^1} ||u_{\eps}(\tau,\cdot)||_{L^p} d\tau.$$
Setting the regularization we obtain
$$
||\partial_x^{\beta}u_{\eps}(t,\cdot)||_{L^p} \leq
\left\{\begin{array}{ll} C, \; \; \; \beta=0, \; n=1\\
C(ln|ln\eps|)^m , \; m>n/2+\beta-1, \; \beta \geq 0, \; n\geq
2\end{array}\right. |ln\eps|^{n(a-1+\theta/2)} $$ $$+ \int_0^t
\left\{\begin{array}{ll} C, \; \; \; \beta=0, \; n=1\\
C(ln|ln\eps|)^m, \; m>n/2+\beta-1, \; \beta \geq 0, \; n\geq 2
\end{array}\right.  |ln\eps|^{n(c-1)}
||u_{\eps}(\tau,\cdot)||_{L^p} d\tau.$$ The moderateness of
$u_{\eps}(t,x)$ from the first step yields

$$||\partial_x^{\beta}u_{\eps}(t,\cdot)||_{L^p} \leq
\left\{\begin{array}{ll} C, \; \; \; \beta=0, \; n=1\\
C(ln|ln\eps|)^m , \; m>n/2+\beta-1, \; \beta \geq 0, \; n\geq
2\end{array}\right. |ln\eps|^{n(a-1+\theta/2)} $$
$$+(CT\left\{\begin{array}{ll} C, \; \; \; \beta=0, \; n=1\\
C(ln|ln\eps|)^m, \; m>n/2+\beta-1, \; \beta \geq 0, \; n\geq 2
\end{array}\right.  |ln\eps|^{n(c-1)}\eps^{-N})
\leq C \eps^{-N}$$ $\exists N>0, \; x\in {\bf R}^n, \; t\in [0,T),
\; 1\leq \theta \leq 1,$ $c<1+1/n.$
 Thus, $u_{\eps}\in {\cal E}_{C^1,
L^p}([0,T), {\bf R}^n),$ $1\leq p \leq \infty.$

Concerning the uniqueness we should to solve the equation
$$
\frac{1}{i}\partial_tW_{\eps}(t,x)=(\triangle-V_{\eps}(x))
W_{\eps}(t,x)+N_{\eps}(t,x)$$
$$W_{\eps}(0,x)=N_{0\eps}(x)\in {\cal N}_{p,q}({\bf R}^n), \;
N_{\eps}(t,x) \in {\cal N}_{C^1, (L^p,L^q)}([0,T), {\bf R}^n), \;
1\leq p,q \leq \infty,$$
$W_{\eps}(t,x)=u_{1\eps}(t,x)-u_{2\eps}(t,x)$ where
$u_{1\eps}(t,x)$ and $u_{2\eps}(t,x)$ are two regularized
solutions to the equation (\ref{donau}). In integral form we have
$$
W_{\eps}(t,x)=S_{n\eps}(t,x)*N_{0\eps}(x)+\int_0^t \int_{{\bf
R}^n} S_{n\eps}(t-\tau, x-y) V_{\eps}(y) W_{\eps}(\tau,y)dy
d\tau$$ $$+\int_0^t \int_{{\bf R}^n} S_{n\eps}(t-\tau, x-y)
N_{\eps}(\tau,y)dy d\tau.$$ Then,
$$
||W_{\eps}(t,\cdot)||_{L^q} \leq
||S_{n\eps}(t,\cdot)||_{L^{\infty}}
||N_{0\eps}(\cdot)||_{L^{q/(q+1)}}+\int_0^t ||S_{n\eps}(t-\tau,
x-\cdot)||_{L^{\infty}} ||V_{\eps}(\cdot)||_{L^1}$$
$$
||W_{\eps}(\tau,\cdot)||_{L^q} d\tau+\int_0^t ||S_{n\eps}(t-\tau,
x-\cdot)||_{L^{\infty}} ||N_{\eps}(\tau,\cdot)||_{L^{q/(q+1)}}
d\tau.$$ Setting the regularization to obtain
$$
||W_{\eps}(t,\cdot)||_{L^q} \leq \left\{\begin{array}{ll} C, \; \;
\;
 n<2\\
C(ln|ln\eps|)^m , \; m> n/2-1, \; \; n \geq 2\end{array}\right.
\eps^s$$ $$+\int_0^t \left\{\begin{array}{ll} C, \; \; \;
 n<2\\
C(ln|ln\eps|)^m , \; m> n/2-1, \; \; n \geq 2\end{array}\right.
|ln\eps|^{n(c-1)} ||W_{\eps}(\tau,\cdot)||_{L^q} d\tau.$$ Gronwall
inequality yields
$$
||W_{\eps}(t,\cdot)||_{L^q} \leq \left\{\begin{array}{ll} C, \; \;
\;
 n<2\\
C(ln|ln\eps|)^m , \; m> n/2-1, \;  \; n \geq 2\end{array}\right.
\eps^s$$ $$ \exp{(CT \left\{\begin{array}{ll} C, \; \; \;
 n<2\\
C(ln|ln\eps|)^m , \; m> n/2-1, \;  \; n \geq 2\end{array}\right.
|ln\eps|^{n(c-1)})}\leq C\eps^s$$ $\forall s \in {\bf N}, \; x\in
{\bf R}^n, \; t\in [0,T), \; \eps<\eps_0,$ $c<1+1/n.$

For the $\beta^{th}$-derivative, $\beta \in {\bf N}_0^n,$ $\beta
\geq 1,$
 we have
$$
\partial_x^{\beta}u_{\eps}(t,x)=
\partial_x^{\beta}S_{n\eps}(t,x)*N_{0\eps}(x)+\int_0^t \int_{{\bf
R}^n}
\partial_x^{\beta}S_{n\eps}(t-\tau,x-y) V_{\eps}(y)
u_{\eps}(\tau,y)dy d\tau$$ $$+ \int_0^t \int_{{\bf R}^n}
\partial_x^{\beta} S_{n\eps}(t-\tau, x-y)  N_{\eps}(\tau,y) dy
d\tau,$$ and
$$
||\partial_x^{\beta}W_{\eps}(t,\cdot)||_{L^q} \leq
||\partial_x^{\beta}S_{n\eps}(t,\cdot)||_{L^{\infty}}
||N_{0\eps}(\cdot)||_{L^{q/(q+1)}}+ \int_0^t ||\partial_x^{\beta}
S_{n\eps}(t-\tau, x-\cdot)||_{L^{\infty}}$$
$$
||V_{\eps}(\cdot)||_{L^1} ||W_{\eps}(\tau,\cdot)||_{L^q}d\tau +
\int_0^t ||\partial_x^{\beta} S_{n\eps}(t-\tau,
x-\cdot)||_{L^{\infty}} ||N_{\eps}(\tau,\cdot)||_{L^{q/(q+1)}}
d\tau.$$ By (\ref{kale}) we have
$$
||\partial_x^{\beta}W_{\eps}(t,\cdot)||_{L^q} \leq
\left\{\begin{array}{ll} C, \; \;  \; \beta=0, \; n=1\\
C(ln|ln\eps|)^m, \; m> n/2+\beta-1, \; \beta \geq 0, \; n\geq 2
\end{array}\right. \eps^s $$
$$+\int_0^t
\left\{\begin{array}{ll} C, \; \; \; \beta=0, \; n=1\\
C(ln|ln\eps|)^m, \; m> n/2+\beta-1, \; \beta \geq 0, \; n\geq 2
\end{array}\right.
|ln\eps|^{n(c-1)}   ||W_{\eps}(\tau,\cdot)||_{L^q} d\tau.$$ The
first step of the induction  yields
$$
||\partial_x^{\beta}W_{\eps}(t,\cdot)||_{L^q} \leq
\left\{\begin{array}{ll} C, \; \; \; \beta=0, \; n=1\\
C(ln|ln\eps|)^m, \; m> n/2+\beta-1, \; \beta \geq 0, \; n\geq 2
\end{array}\right. \eps^s$$
$$
+(CT
\left\{\begin{array}{ll} C, \; \; \;  \beta=0, \; n=1\\
C(ln|ln\eps|)^m, \; m> n/2+\beta-1, \; \beta \geq 0, \; n\geq 2
\end{array}\right. \eps^s
|ln\eps|^{n(c-1)} ) \leq C \eps^s$$ $\forall s \in {\bf N}, \;
x\in {\bf R}^n, \; t\in [0,T), \; \eps< \eps_0,$ under the
condition $c<1+1/n$ where we supposed the
$(ln|ln\eps|)$-boundedeness of the heat semigroup.
 $\Box$

\begin{rem}
 Another
possibility is to set $(ln|ln\eps|)$-boundedeness of potential
$V_{\eps}(x)$ and $(|ln\eps|)$-boundedeness of Schr\"odinger
semigroup under the appropriate condition on the growth of the
mollifier, or put the both to be $(ln|ln\eps|)$-bounded.
Everything leads to the same conclusion: existence-uniqueness
results in Colombeau vector type spaces ${\cal G}_{C^1, (L^q,
L^q)}([0,T), {\bf R}^n)$ for those $1 \leq p,q \leq \infty$ for
which corresponding space is an algebra with multiplication.
\end{rem}

\end{document}